\documentclass[11pt, a4]{article}

\usepackage{amsmath}
\usepackage{amssymb}
\usepackage{amscd}

\usepackage{url}

\usepackage{endnotes}
\let\footnote=\endnote

\usepackage{graphicx}
\usepackage{float}

\usepackage{bmpsize}

\usepackage{array}

\usepackage{authblk}

\newtheorem{theorem}{Theorem}[section]
\newtheorem{lemma}[theorem]{Lemma}

\newtheorem{corollary}[theorem]{Corollary}

\newtheorem{proposition}[theorem]{Proposition}

\title{\bf On Bergman's Diamond Lemma for Ring Theory}
\author{Takao Inou\'{e}}
\date{}
\affil{AI Nutrition Project\\
Artificial Intelligence Center for Health and Biomedical Research,\\
National Institutes of Biomedical Innovation, \\Health and Nutrition, Osaka, Japan\footnote{t.inoue@nibiohn.go.jp \\ \indent 
(Personal e-mail) takaoapple@gmail.com }}

\begin{document}

\maketitle

\begin{abstract}
  This expository and review paper deals with the Diamond Lemma for ring theory, which is proved in the first 
section of G.~M.~Bergman, The Diamond Lemma for Ring Theory, Advances in Mathematics, 29 (1978),  
pp. 178--218.  No originality of the present note is claimed on the part of the author, except for some 
suggestions and figures.  Throughout this paper, I shall mostly use Bergman's expressions in his paper. In Remarks and Notes, the reader will find some useful information on this topic.
\end{abstract}


%%%%%%%%%%%%%%%%%%%%%%%%%%%%%%%

\section{Introduction}

This is an expository and review paper which deals with the Diamond Lemma for ring theory, which is proved in the first section of G.~M.~Bergman, The Diamond Lemma for Ring Theory, Advances in Mathematics, 29 (1978),  pp. 178--218.  No originality of the present note is claimed on the part of the author, except for some suggestions and figures.  Throughout this paper, I shall mostly use Bergman's expressions in his paper. 
In Remarks and Notes, the reader will find some useful information on this topic.

Suppose that $R$ is an associative algebra with 1 over the commutative ring $k$, and that 
we have a presentation of $R$ by a family $X$ of generators and a family $S$ of relations. 
Suppose that each relation $\sigma \in S$ has been written in the form $W_{\sigma} = f_{\sigma}$, 
where $W_{\sigma}$ is a monomial (a product of elements of $X$) and $f_{\sigma}$ is a $k$-linear 
combination of monomials, and that we want to use these relations as instructions for reducing expressions 
$r$ for elements of $R$. That is, if any of the monomials occurring in the 
expression $r$ contains one of the $W_{\sigma}$ as a subword, we substitute 
$f_{\sigma}$ for that subword, and we iterate this procedure as long as possible. In general, 
this process is not always well defined: at each step we must choose \it which \rm reduction 
to apply to \it which \rm subword of \it which \rm monomial. Etcetera. So we are naturally led 
to the following questions:

(1) Under what conditions will such a procedure bring every expression to a unique irreducible form?

(2) Suppose that we have a set of suitable conditions satisfying (1). Does this yield then a 
canonical form for elements of $R$?

The Diamond Lemma is a general result of this sort due to Newman \cite{newman}, which was 
obtained in a graph-theoretic context.\footnote{As a remark for logicians, Newman's paper \cite{newman} 
is closely related to the theory of $\lambda$-calculi. This article also 
contains an interesting observation about a relation of weak Church-Rosser and Church-Rosser 
properties (see Barendregt \cite[p. 58]{henk}).} Let $G$ be an oriented graph. Here the 
vertices of $G$ may be thought as expressions for the elements of some algebraic object (in 
our case, an associative algebra with 1 over the commutative ring $k$) and the edges as reduction 
steps (in our case, reductions using such a rule as $W_{\sigma} = f_{\sigma}$) going from one 
such expression to another one. Newman's result is the following.\footnote{Newman's original formulation 
and terminology for the Diamond Lemma is different from the one in the introduction of the present note 
(see Newman \cite{newman} for the details.)} suppose that 

(i) The oriented graph satisfies \it the descending chain condition\rm . That is, all positively oriented path in $G$ 
terminate; and 

(ii) Whenever two edges, $e$ and $e'$, proceed from one vertex $a$ of $G$, there exists 
positively oritened paths $p$, $p'$ in $G$ leading from the end points $b$, $b'$ of these 
edges to a common vertex $c$. (This condition is called \it the diamond condition\rm .)

%\ref{fig:figure1}
\begin{figure}[h]
 \centering
 \includegraphics[width=8cm]{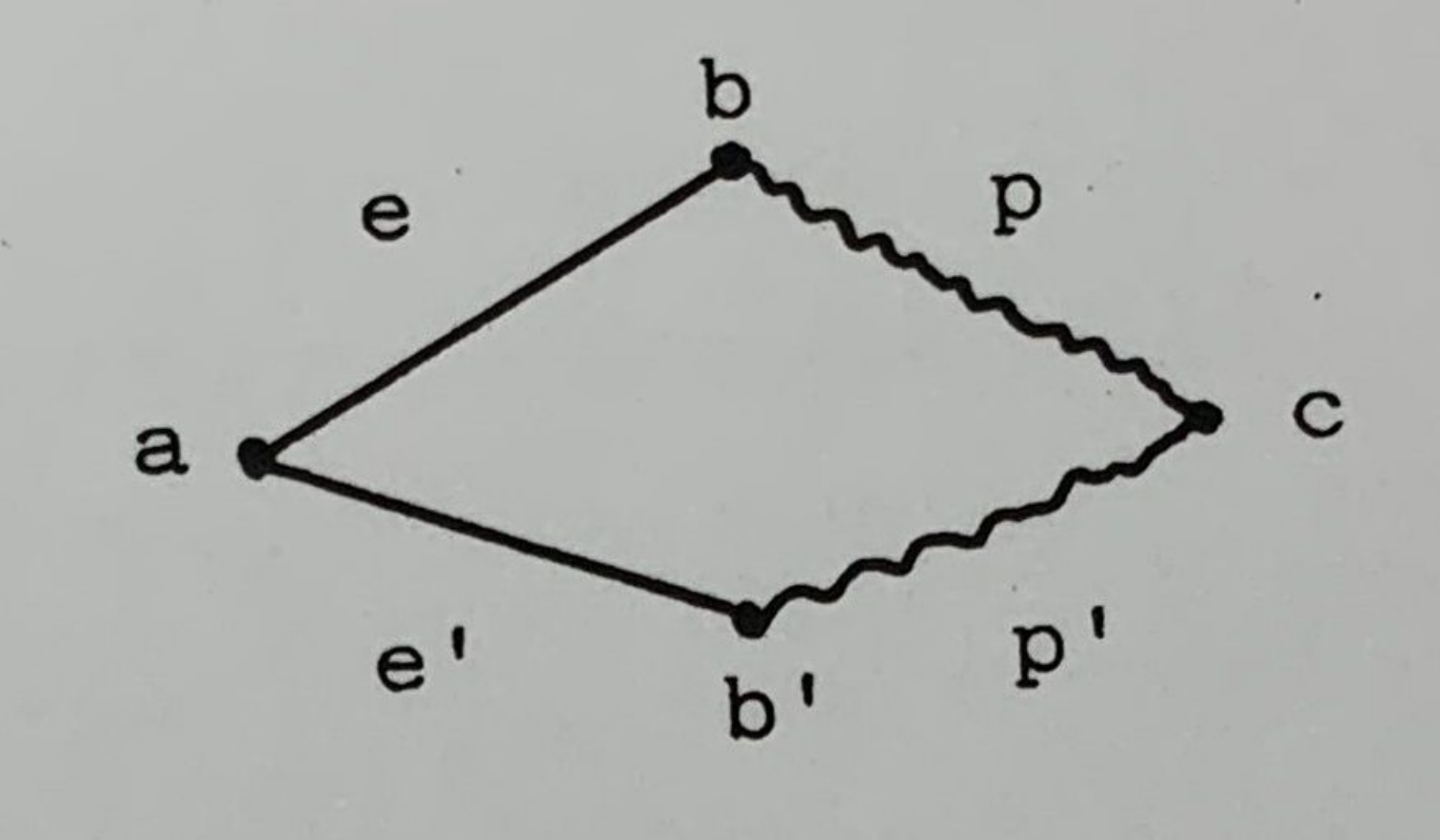}
\caption{Does it look like a diamond?}
 \label{fig:figure1}
\end{figure}%

\it Then \rm every connected component $C$ of $G$  has a unique minimal vertex $m_C$. 

%\ref{fig:figure2}
\begin{figure}[h]
 \centering
 \includegraphics[width=15cm]{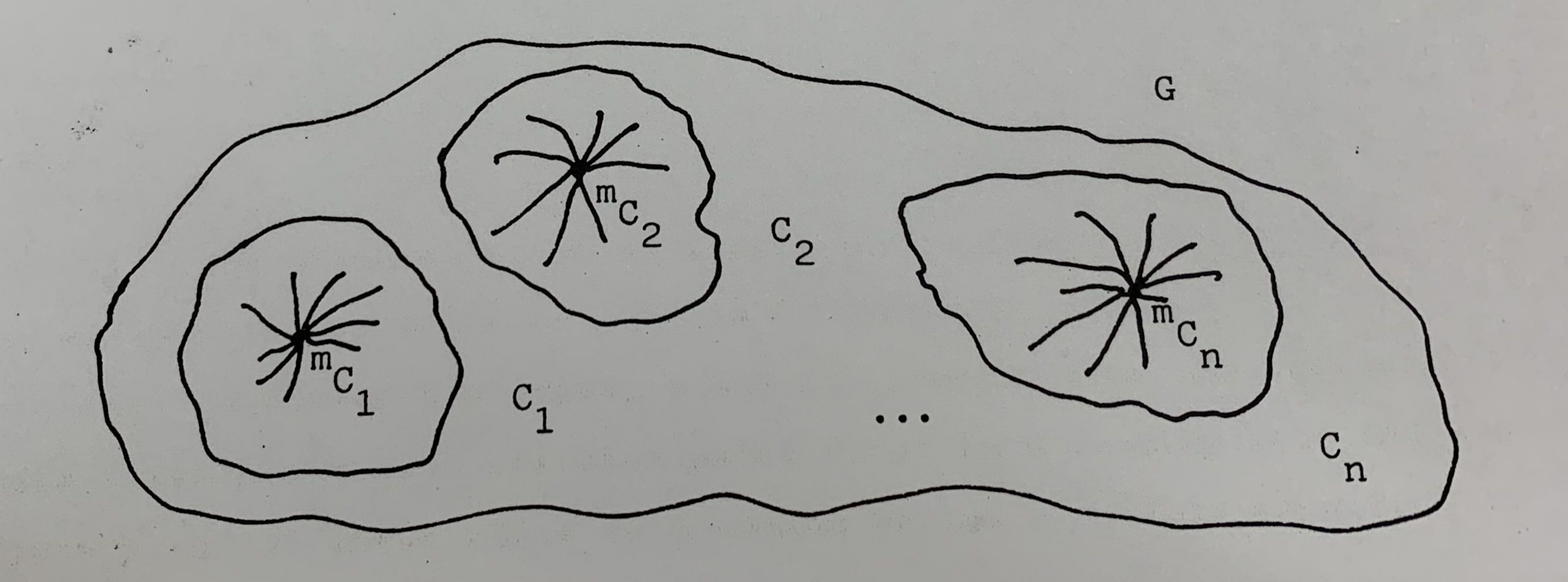}
\caption{The connected components of $G$}
 \label{fig:figure2}
\end{figure}%

\noindent This means that every maximal positively oriented path beginning at a point of $C$ will 
terminate at $m_C$; is other words (in our context) that the given reduction procedures yield unique canonical forms 
for elements of the original algebraic object. 

The main theorem to be proved in the third section, namely the Diamond Lemma for Ring Theory, is an 
analogue of the above observations for the case of associative rings, with reduction procedures of the 
form mentioned earlier. (For our argument in the sequel, we do not follow Newman's graph-theoretic formulation.)

In the following section 2, we introduce a lot of definitions and prove some lemmas and propositions used 
for the proof of the Diamond Lemma. In the last fourth section, we give some suggestions on literatures 
and so on. We have Notes and Appendix at the end of this paper.

%%%%%%%%%%%%%%%%%%%%%%%%%%%%%%%%%%%%%%

\section{Preliminaries}

Let $k$ be a commutative associative ring with $1$, $X$ a set, $<X>$ the free semigroup with $1$ on $X$, and $k<X>$ the
free associative $k$-algebra on $X$, which is the semigroup algebra of $<X>$ over $k$.\footnote{Given such $k$, $X$ 
in the context, the semigroup algebra of 
$<X>$ over $k$ is called the \it free \rm (associative) \it $k$-algebra \rm on $X$.}

Let $S$ be a set of pairs of the form 
$$\sigma = (W_\sigma, f_\sigma)$$ 
\noindent where $W_\sigma \in <X>$, 
$f_\sigma \in k<X>$. For any $\sigma \in S$ and $A, B \in <X>$, let $r_{A\sigma B}$ denote the \it $k$-module \rm 
endomorphism of $k<X>$ that fixes all elements of $<X>$ other than $AW_\sigma B$, and that sends this basis 
element to $Af_\sigma B$. We call the given set $S$ a \it reduction system\rm , and the maps $r_{A\sigma B} : k<X> 
\rightarrow k<X>$ \it reductions\rm .

We say that a reduction $r_{A\sigma B}$ acts \it trivially \rm on an element $a \in k<X>$ if the coefficient of 
$r_{A\sigma B}$ in $a$ is zero. An element $a \in k<X>$ is said to be \it irreducible \rm if every reduction acts 
trivially on $a$.

\begin{proposition} The irreducible elements of $k<X>$ form a $k$-submodule of $k<X>$, denoted by $k<X>_{irr}$.
\end{proposition}
\sc Proof\rm . Let $a$, $b$ be any irreducible elements of $k<X>$ and $\lambda$ any element of $k$. Let $r$ 
be a reduction, say $r = r_{A\sigma B}$. The coefficient of $r_{A\sigma B}$ in $a$ and $b$ is zero, respectively. 
Thus so is that of $r_{A\sigma B}$ in $a - b$ and $\lambda a$. Trivially $0$ is irreducible. This completes the 
proof. $\Box$

A finite sequence of reduction $r_1, \dots , r_n (r_i= r_{A_i{\sigma_i}B_i})$ is said to be \it final \rm on $a \in k<X>$ 
if $r_n \cdots r_1(a) \in k<X>_{irr}$.

An element $a$  of $k<X>$ is called \it reduction-finite \rm if for every infinite sequence $r_1, r_2, \dots$ of reduction, 
$r_i$ acts trivially on $r_{i - 1} \cdots r_1(a)$ for all sufficiently large $i$. If $a$ is reduction-finite, then 
any maximal sequence of reductions $r_i$ acts \it nontrivially \rm on $r_{i - 1} \cdots r_1(a)$ is finite, and hence 
a final sequence. 

\begin{proposition} The reduction-finite elements of $k<X>$ form a $k$-submodule of $k<X>$.
\end{proposition}
\sc Proof\rm .  Suppose that $a$ and $b$ are reduction-finite elements and $\lambda$ an element of $k$. 
Then there are natural number $i$ and $i$ such that for every infinite sequence $r_1, r_2, \dots$ of reductions, 
$r_i$ and $r_j$ act trivially on $r_{i - 1} \cdots r_1(a)$ and $r_{j - 1} \cdots r_1(b)$, respectively. Take $l = max(i, j)$. 
For every infinite sequence $r_1, r_2, \dots$ of reductions, $r_l$ and $r_i$ act trivially on 
 $r_{l - 1} \cdots r_1(a - b)$ and $r_{j - 1} \cdots r_1(\lambda a)$, respectively. Thus, $a - b$ and $\lambda a$ are
reduction-finite. $0$ is clearly reduction-finite. This completes the proof. $\Box$

We call an element $a \in k<X>$ \it reduction-unique \rm if 

(1) it is reduction-finite; and 

(2) its images under all final sequences are the same. (This common value is denoted by $r_s(a)$.)

\begin{lemma} 

(i) The set of reduction-finite elements of $k<X>$ form a $k$-submodule of $k<X>$, and $r_s$ is a $k$-linear 
map of this submodule into $k<X>_{irr}$.

(ii) Suppose $a, b, c \in k<X>$ are such that for all monomials $A, B, C$ occurring with nonzero coefficient in $a, b, c$, 
respectiely, the product $ABC$ is reduction-unique. (In particular this implies that $abc$ is reduction-unique.) 
Let $r$ be any finite composition of reductions. Then $ar(b)c$ is reduction-unique, and $r_s(ar(b)c) = r_s(abc)$. 
(Note that 'finite' means 0 or $\geq$ 2. When $r$ is a single reduction, $ar(b)c$ should have the same property 
as that of $abc$.)

\end{lemma}
\sc Proofr\rm .  (i) Suppose that $a, b \in k<X>$ are reduction-unique, and $\alpha \in k$. By Proposition 2.2, 
$\alpha a + b$ is reduction-finite. Let $r$ be any composition of (finite) reductions final on $\alpha a + b$. 
Since $a$ is reduction-unique, we can find a composition of (finite) reduction $r'$ such that $r'r(a) = r_s(a)$, 
and similarly there is a composition of reductions $r''$ such that  $r''r'r(b) = r_s(b)$. Because 
$r(\alpha a + b) \in k<X>_{irr}$, we have 
\begin{equation*}
\begin{split}
r(\alpha a + b) &= r''r'r(\alpha a + b)\\
&= \alpha r''r'r(a) + r''r'r(b)\\
&= \alpha r''r_s(a) + r_s(b)\\ 
&= \alpha r_s(a) + r_s(b)
\end{split}
\end{equation*}
That is, images of $\alpha a + b$ under all such final sequenes of reductions are the same, i.e. $\alpha r_s(a) + r_s(b)$. 
Thus, $\alpha a +b$ is reduction-unique and so is $b - a$ with $\alpha = -1$. $0$ is clearly reduction-unique. 
Therefore, the set of reduction-unique elements of $k<X>$ forms a $k$-submodule of $k<X>$.  Since 
$r_s(\alpha a + b) = r(\alpha a + b)$, $r_s(\alpha a + b) = \alpha  r_s(a) + r_s(b)$. For any reduction-unique elements 
of $k<X>$, $r_s(s) \in k<X>_{irr}$ is 
clear. Thus, $r_s$ is a $k$-linear map of the module into  $k<X>_{irr}$.

(ii) Suppose that the assumption of (ii) holds. And say 
$$a = \sum_i \alpha_i A_i, \enspace  b = \sum_j \beta_j B_i, \enspace c = \sum_l \gamma_l C_l.$$
So $abc = \sum_{i, j, l} \alpha_i \beta_j \gamma_l A_iB_jC_l$ and for any triple $(i, j, l)$, $A_iB_jC_l$ is reduction-unique. 
So $abc$ is reduction-unique and $r_s(abc) =  \sum_{i, j, l} \alpha_i \beta_j \gamma_l r_s(A_iB_jC_l)$. Let $r$ be any finite 
composition of reductions. It is sufficient to consider the case where $a, b, c$ are monomials and $r$ is a single 
reduction $r_{D\sigma E}$. In this case, $Ar_{D\sigma E}(B)C = r_{AD\sigma EC}$, which is the image of $ABC$ under a 
reduction, hence is reduction-unique if $ABC$ is so, with the reduced form. $\Box$

By an \it overlap ambiguity \rm of $S$ we mean a 5-tuple $(\sigma, \tau, A, B, C)$ with $\sigma, \tau \in S$ and 
$A, B, C \in <X> - \{1\}$, such that $W_\sigma = AB$, $W_\tau = BC$. We say that the overlap ambiguity 
$(\sigma, \tau, A, B, C)$ is \it resolvable \rm if there exist compositions of reductions, $r$, $r'$, such that 
$r(f_\sigma C) = r'(Af_\tau)$: in other words, $f_\sigma C$ and $Af_\tau$ can be reduced to a common expression 
(This corresponds to the diamond condition seen in the introduction).

%\ref{fig:figure3}
\begin{figure}[h]
 \centering
 \includegraphics[width=10cm]{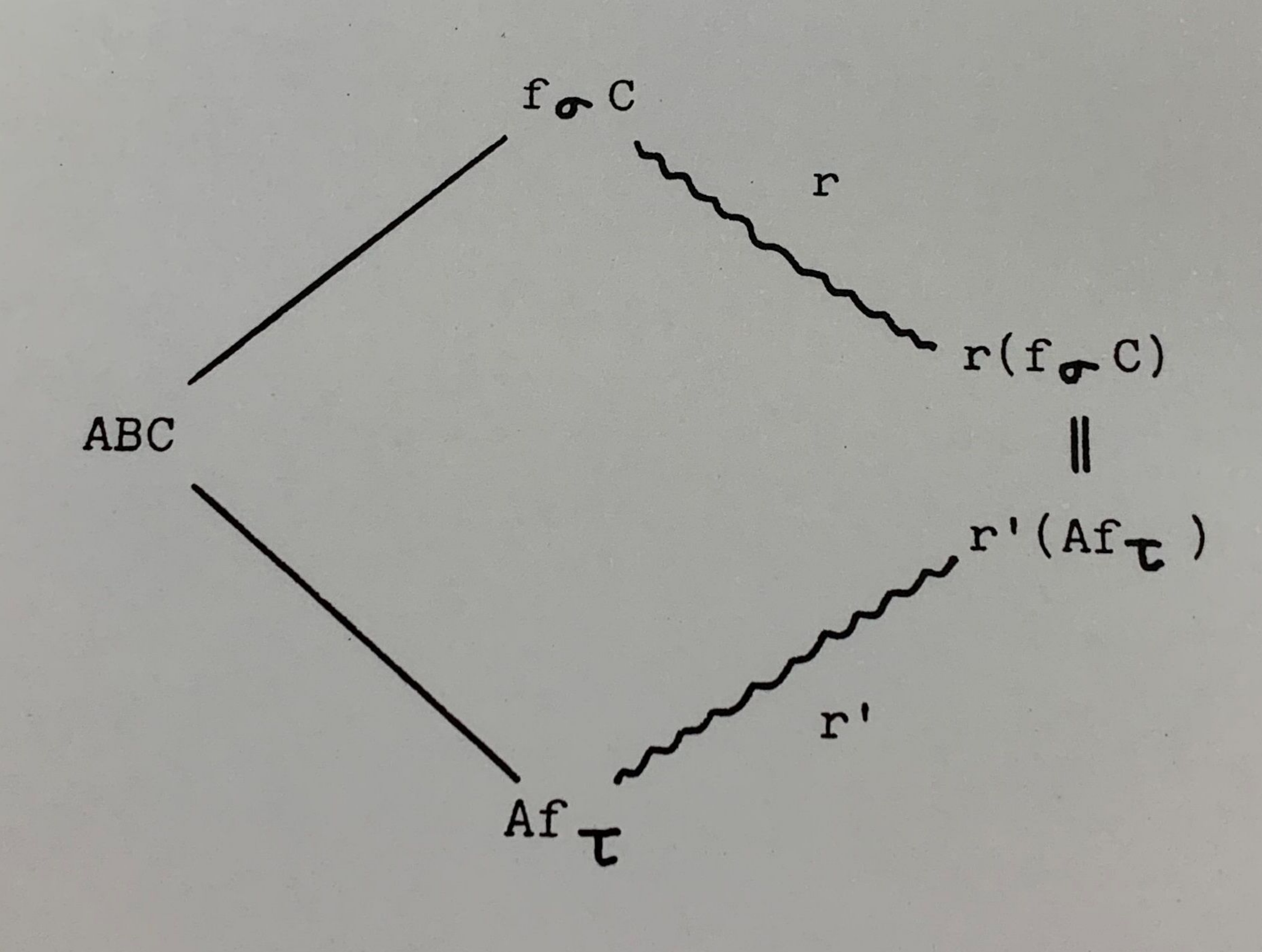}
\caption{The diamond condition?}
 \label{fig:figure3}
\end{figure}%

Similarly, a 5-tuple $(\sigma, \tau, A, B, C)$ with $\sigma \neq \tau \in S$ and $A, B, C \in <X>$ is called 

and \it inclusion ambiguity\rm \footnote{Inclusion ambiguities are, in a sense, always avoidable. Suppose that $S$ is a reduction 
system for a free algebra $k<X>$. Let us construct a subset $S' \subseteq S$ by (1) deleting all $\sigma \in S$ such that 
$W_\sigma$ contains a \it proper \rm subword of the form $W_\tau (\tau \in S)$, and (2) whenever more than one element 
$\sigma_1, \sigma_2, \dots \in S$ act on the same monomial (i.e. $W_{\sigma_1} = W_{\sigma_2} = \dots$) dropping all but one of 
the $\sigma_i$ from $S$. Then $S' \subseteq S$ will have the property that $a \in k<X>$ is reducible under $S$. But from 
this it follows that if $a \in k<X>$ is reduction-unique under $S$, then so is $a$ under $S'$ and $r_S'(a) = r_S(a)$. Hence if $S$ is 
such that every element of $k<X>$ is reduction-unique under it, then $S'$ has the same property and $r_{S'} = r_S$. Therefore, $S'$, 
which has no inclusion ambiguities, defines the same ring and the same canonical form as $S$. (This remark is due to 
Bergman \cite[p. 192]{bergman1}.) } 
if $W_\sigma = B$, $W_\tau = ABC$. The inclusion ambiguity is called \it resolvable \rm if there exists compositions of 
reductions, $r$ and $r'$, such that $r(Af_\sigma B) = r'(f_\tau)$. 

By a \it semigroup partial ordering on \rm $<X>$, we mean a partial order "$<$" such that 
$B < B' \Leftrightarrow ABC <AB'C$ for any $A, B, B', C \in k<X>$, and it is called \it compatible with \rm $S$ if 
for all $\sigma \in S$, $f_\sigma$ is a (finite) linear combination of monomials $< W_\sigma$. 

If $\leq$ is a semigroup partial ordering on $<X>$ compatible with the reduction system $S$, and $A$ is any element 
of $<X>$, let $I_A$ denote the submodule of $k<X>$ spanned by all elements $B(W_\sigma - f_\sigma)C$ such that $BW_\sigma C < A$. 
We say that an overlap  (inclusion) ambiguity  $(\sigma, \tau, A, B, C)$ is \it resolvable relative 
to \rm $\leq$ if $f_\sigma C - Af_\tau \in I_{ABC}$ $(Af_\sigma C - f_\tau \in I_{ABC})$. 

The following lemma is trivial. But it will be useful in what follows.

\begin{lemma} 
Let $a \in k<X>$. Suppose that $a$ contains a monmial of the form $AW_\sigma B$ with a coefficient 
$\lambda ( \neq 0) \in k$. Then we have $$r_{A\sigma B}(a) = a - \lambda A(W_\sigma - f_\sigma) B.$$
\end{lemma}
\sc Proof\rm . The lemma immediately follows from the following observation:
$$r_{A\sigma B}(\lambda AW_\sigma B) = \lambda Af_\sigma B =  \lambda AW_\sigma B - 
\lambda A(W_\sigma - f_\sigma) B.$$
$\Box$

Let $I$ denote the two-sided ideal of $k<X>$ generated by th elements $W_\sigma - f_\sigma (\sigma \in S)$. 
As a $k$-module, $I$ is spanned by the products $A(W_\sigma - f_\sigma)B$. 

\begin{proposition}
Let $a \in k<X>$. Suppose that $a$ is reduction-unique. Then, if $r_s(a) = 0$, then $a$ is an element of $I$.
\end{proposition}
\sc Proof\rm . If $a = 0$, the $0 \in I$ is trivial. So assume $a \neq 0$. Suppose that $r$ is the composition of 
a sequence final on $a$ and say,
$$r = r_{A_n{\sigma_n}B_n} \cdots r_{A_2{\sigma_2}B_2}r_{A_1{\sigma_1}B_1}.$$
Then, By Lemma 2.4, $r(a)$ is of the form $a - \sum \lambda_i A_i(W_{\sigma_i} - f_{\sigma_i})B_i$ with 
$\lambda_i (\neq 0) \in k$ for all $i$. Because $a$ is reduction-unique, $r(a) = r_s(a)$, which implies 
$a =  \sum \lambda_i A_i(W_{\sigma_i} - f_{\sigma_i})B_i$. Thus, $a \in I$. $\Box$

%%%%%%%%%%%%%%%%%%%%%%%%%%%%%%%%%%%%

\section{The Diamond Lemma}

The following theorem is called the Diamond Lemma for Ring Theory.\footnote{As an application of the lemma, Bergman 
shows, for example an alternative proof of Poincare\'{e}-Birkoff-Witt Theorem in \cite[p.186]{bergman1}, which gives a 
basis of the universal enveloping algebra $U(g)$ of $g$ if we know a basis of a Lie algebra $g$. Also see Varadarajan \cite{Varadarajan}. }

\begin{theorem}
Let $S$ be a reduction system for a free associative algebra $k<X>$ (a subset of $<X> \times k<X>$), and $\leq$ a 
semigroup partial ordering on $<X>$ , compatible with $S$, and satisfies the descending chain condition. Then the 
following conditionsa are equivalent:

(a) All ambiguities of $S$ are resolvable.

(a') All ambiguities of $S$ are resolvable relative to $\leq$.

(b)  All elements of $k<X>$ are reduction-unique under $S$.

(c) A set of representatives in $k<X>$ for the elements of the algebra $R = k<X>/I$ determined by the generators 
and the relations $W_\sigma =f_\sigma$  ($\sigma \in S$) is given by the $k$-submodule $k<X>_{irr}$ sapanned by 
the $S$-irreducible monomials of $<X>$.

When these conditions hold, $R$ may be identified with the $k$-module $k<X>_{irr}$, made a $k$-algebra by 
the multiplication $a \cdot b = r_S(ab)$.
\end{theorem}

\noindent \sc Proof\rm . First we see from our general hypothesis, that every element of $<X>$ is reduction-finite. 
We can prove this formally by induction with respect to the partial ordering with the descending chain condition $\leq$. 
But here we prove it informally to make the situation clearer. For illustrative puposes, suppose that $a \in <X>$ has a 
monomial of the form $AW_\sigma B$ ($\sigma \in S$). By a reduction $r_{A\sigma B}$, $r_{A\sigma B}(a)$ has a monomial  
of the form $Af_\sigma B$. Since $\leq$ is compatible with $S$, $f_\sigma$ is of the form $\sum \lambda_i W_i^\sigma$ 
with $W_i^\sigma < W_\sigma$ for any $i$. So every monomial of $Af_\sigma B$ is of the form $AW_i^\sigma B$ for some $i$. 
If the monomial $AW_i^\sigma B$ contains a subword $W_\tau$ ($\tau \in S$), say $AW_i^\sigma B  = A'W_\tau B'$, by a 
reduction $r_{A'\tau B'}$, $r_{A'\tau B'} r_{A\tau B}(a)$ has a monomial of the form $A'f_\tau B'$. By compatibility of 
$\leq$, $f_\tau$ is again a (finite) linear combination of monomials $< W_i^\sigma$, say $f_\tau = \sum \mu_j W_j^\tau$ 
with $W_j^\tau < W_i^\sigma$ for any $j$. So $r_{A'\tau B'} r_{A\sigma B}(a)$ has a monomial of the form $A'W_j^\tau B'$ 
for all $j$. If we iterate this process, we will get a sequence of monomials, for example 
$W_\sigma > W_i^\sigma > W_j^\tau > \cdots$. All of such  sequences must be finite because of the descending chain 
condition. It is also clear that the number of all the possible sequences is finite. (See the following figure.) 

%\ref{fig:figure4}
\begin{figure}[h]
 \centering
 \includegraphics[width=15cm]{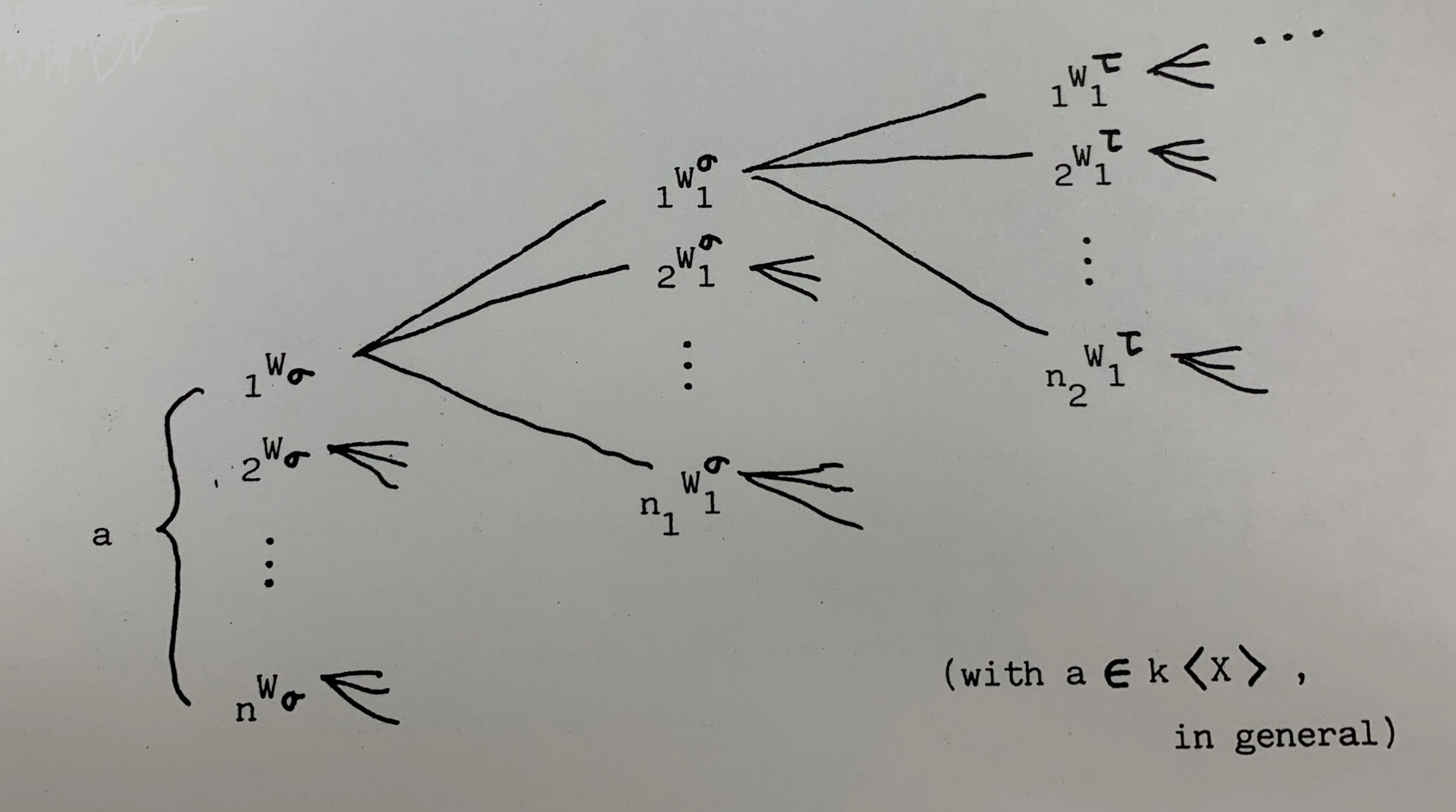}
\caption{The number of all the possible sequences is finite.}
 \label{fig:figure4}
\end{figure}%

\noindent Therefore, $a$ is reduction-finite. Since every element of $<X>$ is reduction-finite, hence so is every element of $k<X>$.

Next we prove (b) $\Leftrightarrow$ (c). We note first that (c) simply says 
$$k<X> \enspace \simeq  \enspace k<X>_{irr} \oplus \enspace I.$$ 
Assuming (b), we show that the following sequence, 
$$(*) \enspace\enspace 0 \longrightarrow I \stackrel{i}{\longrightarrow} k<X> \stackrel{r_S}{\longrightarrow} k<X>_{irr} \longrightarrow 0,$$
is a short exact sequence $k$-module homomorphisms, where $i$ is inclusion map. If so, it is immediate to conclude 
$k<X> \enspace \simeq  \enspace k<X>_{irr} \oplus \enspace I$. That is, it follows from $r_Si = id_{k<X>_{irr}}$ (see Appendix).

Case 1: inclusion map $i$ is injective. So $0 \rightarrow I \rightarrow k<X>$ is exact.

Case 2: $I = im(i) \subseteq ker(r_S)$ is easily seen because for all $A, B, \sigma$, 
\begin{equation*}
\begin{split}
r_S(A(W_\sigma - f_\sigma)B) &= r_S(AW_\sigma B) - r_S(Af_\sigma B)\\
&(\mbox{$r_S$ is $k$-linear. (Lemma 2.3.(i))}\\
&= r_S(AW_\sigma B) - r_S(Ar_{1\sigma 1}(1W_\sigma1)B) = 0\\ 
&\mbox{(By Lemma 2.3.(ii), $A, W_\sigma, B$ are reduction-unique by (b).)}
\end{split}
\end{equation*}
By Proposition 2.5, $ker(r_S) \subseteq im(i) = I$ is obvious. Thus, $im(i) = ker(r_s)$. 
This is nothing but the exactness of $I \rightarrow k<X> \rightarrow k<X>_{irr}$.

Case 3: $r_S$ is surjective, since for any $x \in k<X>_{irr}$, $r_S(x) = x$. This means that $k<X> \rightarrow k<X>_{irr} \rightarrow 0$ is exact.

Conversely, we assume (c) and suppose $a \in k<X>$ can be reduced to either of $b, b' \in K<X>_{irr}$. Then we have 
$b -b' \in k<X>_{irr} \cap I = \{0\}$, which proves (b).

The final comment in the statement of the theorem is now clear. From the above, $R = k<X>/I \simeq k<X>_{irr}$ is immediate 
and $k<X>_{irr}$ is a $k$-algebra with the multiplication $a \cdot b = r_S(ab)$ by the following (1) and (2):

(1) $k<X>_{irr}$ is a $k$-submodule of $k<X>$ by Lemma 2.3.(i). Thus it is a $k$-module.

(2) For any $a, b \in k<X>_{irr}$ and any $\alpha \in k$, 
$$ \alpha (a \cdot b) = (\alpha a) \cdot b = a \cdot (\alpha b),$$
holds, because we see that 
\begin{equation*}
\begin{split}
 (\alpha a) \cdot b &= r_S((\alpha a) \cdot b)\\
&= r_S(\alpha (ab))\\
&\mbox{$((\alpha a)b = \alpha (ab)$ holds since $k<X>$ is a $k$-algebra.)}\\
&= \alpha r_S(ab) \enspace \mbox{($r_S$ is $k$-linear.)}\\
&= \alpha (ab)
\end{split}
\end{equation*} 
and with similar remarks, $a \cdot (\alpha b) = \alpha (a \cdot b)$.

We next deal with the proof of (b) $\Leftrightarrow$ (a). Suppose (b). We consider only the case of overlap ambiguities, 
because those of inclusion ones are similarly taken case of. Let $(\sigma, \tau, A, B, C)$ be any overlap ambiguity. 
$f_\sigma C$, $Af_\tau$ are reduction-unique by (b). So we may take compositions of reductions $r$ and $r'$ which are final on 
$f_\sigma C$ and $Af_\tau$, respectively. By (b), $ABC$ is reduction-unique. Moreover, $rr_{1\sigma C}$ and $r'r_{A\tau 1}$ are 
obviously final on $ABC$. So we see that 
$$r(f_\sigma C) = rr_{1\sigma C}(ABC) = r_S(ABC) = r'r_{A\tau 1}(ABC) = r(Af_\tau).$$
This means that the ambiguity is resolvable.

In this paragraph, we porve (a) $\Rightarrow$ (a').  We assume (a). First we consider the case of overlap ambiguities. 

Let $(\sigma, \tau, A, B, C)$ be any overlap ambiguity of $S$. By (a), it is resolvable. That is, there are compositions of 
reductions $r$ and $r'$ such that $r(f_\sigma C) = r'(Af_\tau)$, say 
$$r = r_{D_n\sigma_nE_n} \cdots r_{D_1\sigma_1E_1} \enspace \mbox{and} \enspace r' = r_{D'_m\tau_mE'_m} \cdots r_{D'_1\tau_1E'_1}.$$
By Lemma 2.4, we obtain, 
$$r(f_\sigma C) = f_\sigma C - \sum_{i = 1}^{n} \lambda_iD_i(W_{\sigma_i} - f_{\sigma_i})E_i$$ 
and
$$r'(Af_\tau) = Af_\tau - \sum_{i = 1}^{m} \mu_iD'_i(W_{\tau_i} - f_{\tau_i})E'_i$$ 
with $\lambda_i (\neq 0), \mu_j (\neq 0) \in k$ for $1 \leq i \leq n$ and $1 \leq j \leq m$. Under the diamond condition 
in our sense, we may have then, 
$$f_\sigma C - Af_\tau =  \sum_{i = 1}^{n} \lambda_iD_i(W_{\sigma_i} - f_{\sigma_i})E_i - \sum_{i = 1}^{m} \mu_iD'_i(W_{\tau_i} - f_{\tau_i})E'_i.$$
Further, it is not so difficult to (I) $D_iW_{\sigma_i}E_i < W_\sigma C = ABC$ and (II) $D'_jW_{\tau_j}E'_j < W_\tau = ABC$ for any 
$1 \leq i \leq n$ and any $1 \leq j \leq m$. For the verification of (I), we show only the case of $i = 1$. The rest of the proof is taken care 
of tby induction. Since $\leq$ is compatible with $S$, if $f_\sigma$ is of the form $\sum \alpha_iZ_i$, then 
$Z\i < W_\sigma$ holds for any $i$. Further, $Z_i < W_\sigma$ leads to $Z_iC < W_\sigma C = ABC$ for all $i$. So we must have 
$D_1W_{\sigma_1}E_i = Z_iC$ for some $i$. Hence, $D_1W_{\sigma_1}E_1 < ABC$. We can verify (II) similarly. So we omit the verification. 
For inclusion ambiguity, we can also show resolvability relative to $\leq$ in a completely similar way. So this is left to the reader.

In this paragraph, we take care of the last implication to be shown, i.e. (a') $\Leftarrow$ (b). It suffices to prove all monomials 
$D \in <X>$ reduction-unique, since the reduction-unique elements of $k<X>$ form a submodule (Lemma 2.3.(i)). That is, if every 
monomial of $<X>$ is reduction-unique, then $k<X>_{irr} = k<X>$. We assume inductively that 
all monomials $< D$ are reduction-unique. Thus the domain of $r_S$ includes the submodule spanned by all these monomials, 
so the kernel of $r_S$ contains $I_D$. That is, if $a \in ker(r_S)$, then by Proposition 2.5, $a$ is of the form 
$\sum \lambda_iA_i(W_{\sigma_i} - f_{\sigma_i})B_i$ with $A_iW_{\sigma_i}B_i < D$ for any $i$, which means $a \in I_D$. We must now 
show that given any two reductions $r_{L\sigma M'}$ and $r_{L'\tau M}$ each acting nontrivially on $D$ (and hence each sending $D$ to a 
linear combination of monomials $< D$), we will have 
$$r_S(r_{L\sigma M'}(D)) = r_S(r_{L'\tau M}(D)).$$ 
We have to check three case for that, according to the relative locations of the subwords $W_\sigma$ and $W_\tau$ in the monomial $D$. 
We may assume without loss of generality that $length(L) \leq length(L')$, in other words, that the indicated copy of $W_\sigma$ in $D$ 
begins no later than that of $W_\tau$.

Case 1: The subwords $W_\sigma$ and $W_\tau$ overlap in $D$, neither contains the other, figured as follows: under the condition 
$length(L) \leq length(L')$, 

\begin{center}
\begin{tabular}{|>{\centering}p{2cm}|>{\centering}p{2cm}|>{\centering}p{1cm}|} \hline 
$L$ & $W_\sigma$ & $M'$  
\end{tabular}
\end{center}

\begin{center}
\begin{tabular}{|>{\centering}p{2.5cm}|>{\centering}p{2cm}|>{\centering}p{0.5cm}|} \hline 
$L'$ & $W_\tau$ & $M$  
\end{tabular}
\end{center}

\noindent Tnen $D = LABCM$, where $(\sigma, \tau, A, B, C)$ is an overlap ambiguity of $S$, i.e. 
$W_\sigma = AB, W_\tau = BC, \sigma, \tau \in S, A, B, C \in <X> - \{1\}$. Then,  
\begin{equation*}
\begin{split}
F &:= r_{L\sigma M'}(D) - r_{L'\tau M}(D)\\
&= Lf_\sigma CM -LAf_\tau M\\
&= L(f_\sigma CM -Af_\tau M)\\
&= L(f_\sigma C -Af_\tau)M. \enspace \enspace \cdots \cdots \enspace \enspace (1.1)
\end{split}
\end{equation*} 
By (a') every overlap ambiguity is resolvable relative to $\geq$. So we have $f_\sigma C - Af_\tau \in I_{ABC}$ by definition. 
That is, 
$$f_\sigma C - Af_\tau = \sum \lambda_iD_i(W_{\sigma_i} - f_{\sigma_i})E_i$$
with $D_iW_{\sigma_i}E_i < ABC$ for any $i$. Substitute this to (1.1). Then we get,  
$$F = \sum \lambda_iLD_i(W_{\sigma_i} - f_{\sigma_i})E_iM. \enspace \enspace \cdots \cdots \enspace \enspace (1.2)$$
Since $\geq$ is a semigroup ordering, the following inequality, 
$$LD_i(W_{\sigma_i})E_iM < LABCM. \enspace \enspace \cdots \cdots \enspace \enspace (1.3)$$
holds by $D_i(W_{\sigma_i})E_i < ABC$. From (1.2) and (1.3), it follows that $F \in I_{LABCM}$. Thus $r_S(F) = 0$, in other words,
$$r_S(r_{L\sigma M'}(D) - r_{L'\tau M}(D)) = 0,$$
so $r_S(r_{L\sigma M'}(D)) = r_S(r_{L'\tau M}(D))$.

The next case is similarly dealt with as the case 1. But we shall work it out for the sake of the reader.

Case 2: One of the subwords $W_\sigma$, $W_\tau$ is contained in the other. By 
$length(L) \leq length(L')$, we have the following case where $W_\sigma$ contains $W_\tau$, figured below.

\begin{center}
\begin{tabular}{|>{\centering}p{1cm}|>{\centering}p{2.5cm}|>{\centering}p{1.5cm}|} \hline 
$L$ & $W_\sigma$ & $M'$  
\end{tabular}
\end{center}

\begin{center}
\begin{tabular}{|>{\centering}p{2cm}|>{\centering}p{1cm}|>{\centering}p{2cm}|} \hline 
$L'$ & $W_\tau$ & $M$  
\end{tabular}
\end{center}

$$$$

\noindent Then $D = LABCM'$, $CM' = M$ and $L' = LA$, where $(\sigma, \tau, A, B, C)$ is an inclusion ambiguity of $S$, i.e. 
$W_\tau = B, W_\sigma = ABC$ with $\tau \neq \sigma \in S, A, B, C \in <X>$. Then,  
\begin{equation*}
\begin{split}
F &:= r_{L\sigma M'}(D) - r_{L'\tau M}(D)\\
&= Lf_\sigma M'-LAf_\tau CM'\\
&= L(f_\sigma -Af_\tau C)M'. \enspace \enspace \cdots \cdots \enspace \enspace (2.1)
\end{split}
\end{equation*} 
By (a') we know that every inclusion ambiguity is resolvable relative to $\geq$. So we get $Af_\tau C - f_\sigma \in I_{ABC}$ by definition. 
That is, 
$$Af_\tau C - f_\sigma = \sum \lambda_iD_i(W_{\sigma_i} - f_{\sigma_i})E_i$$
with $D_iW_{\sigma_i}E_i < ABC$ for any $i$. By substituing this to (2.1), we obtain, 
$$F = \sum (-\lambda_i)LD_i(W_{\sigma_i} - f_{\sigma_i})E_iM'. \enspace \enspace \cdots \cdots \enspace \enspace (2.2)$$
Because $\geq$ is a semigroup ordering, the following inequality, 
$$LD_i(W_{\sigma_i})E_iM' < LABCM. \enspace \enspace \cdots \cdots \enspace \enspace (2.3)$$
holds by $D_i(W_{\sigma_i})E_i < ABC$. From (2.2) and (2.3), we get $F \in I_{LABCM'} = I_D$, so  $r_S(F) = 0$. That is,
$$r_S(r_{L\sigma M'}(D) - r_{L'\tau M}(D)) = 0,$$
thus $$r_S(r_{L\sigma M'}(D)) = r_S(r_{L'\tau M}(D)).$$

The following is our last case to check, with which we complet the whole proof of the Diamond Lemma. 

Case 3: We consider the case where $W_\sigma$ and $W_\tau$ is disjoint. By the condition on thelength of $L$ and $L'$, 
$length(L) \leq length(L')$, the case is figured below.

$$$$

\begin{center}
\begin{tabular}{|>{\centering}p{0.7cm}|>{\centering}p{0.6cm}|>{\centering}p{3.7cm}|} \hline 
$L$ & $W_\sigma$ & $M'$  
\end{tabular}
\end{center}

\begin{center}
\begin{tabular}{|>{\centering}p{2.5cm}|>{\centering}p{1.7cm}|>{\centering}p{0.8cm}|} \hline 
$L'$ & $W_\tau$ & $M$  
\end{tabular}
\end{center}

$$$$

\noindent So we may assume $D = LW_\sigma NW_\tau M$, i.e., 

$$$$

$\hspace*{4.4cm} \overbrace{\hspace*{4.8cm}}^{M'}$
\begin{center}
\begin{tabular}{|>{\centering}p{0.7cm}|>{\centering}p{0.6cm}|>{\centering}p{1.2cm}|>{\centering}p{1.7cm}|>{\centering}p{0.8cm}|} \hline 
$L$ & $W_\sigma$ & $N$ & $W_\tau$ & $M$  
\end{tabular}
\end{center}
$\hspace*{2.8cm} \underbrace{\hspace*{3.6cm}}^{L'}$

\bigskip

\noindent is our present case with 
$$r_{L\sigma M'}(D) = Lf_\sigma NW_\tau M, \enspace \enspace (M' = NW_\tau M)$$
and 
$$r_{L'\tau M}(D) = LW_\sigma Nf_\tau M. \enspace \enspace (L' = LW_\tau N)$$
By the general assumption, the ordering $\leq$ is compatible with $S$. So $f_\sigma$ can be witten as a linear combination of 
monomials $<W_\sigma$, say $f_\sigma = \sum \lambda_iZ_i$ with $Z_i < W|_\sigma$, $\lambda_i (\neq 0) \in k$. The ordering 
is a semigroup ordering, so for any $i$, we have $LZ_iNW_\tau M < LW_\sigma NW_\tau M = D$ from $Z_i < W_\sigma$. 
By induction hypothesis, $LZ_iNW_\tau M$ is reduction=unique for all $i$. Let $a = 1, c =1$ and 
$b = Lf_\sigma NW_\tau M = \sum \lambda_iLZ_iNW_\tau M$. Then, for all monomials $A, B, C$ occurring with non-zero coefficient 
in $a, b, c$, i.e. $1, LZ_iNW_\tau M (\forall i), 1$, respectively, the product $ABC$, namely $LZ_iNW_\tau M (\forall i)$ is 
reduction-unique. Apply now Lemma 2.3.(ii) to such $a, b, c$ with $r = r_{(Lf_\sigma N)\tau M}$. Then we have $ar(b)c$ is 
reduction-unique and $r_S(ar(b)c) = r_S(abc)$. In other words, $Lf_\sigma Nf_\tau M$ is reduction-unique and 
$$r_S(Lf_\sigma Nf_\tau M) = r_S(Lf_\sigma NW_\tau M). \enspace \enspace \cdots \cdots \enspace \enspace (3.1)$$
Similarly, we can obtain, 
$$r_S(Lf_\sigma Nf_\tau M) = r_S(LW_\sigma Nf_\tau M). \enspace \enspace \cdots \cdots \enspace \enspace (3.2)$$
From (3.1) and (3.2), it follows that $ r_S(Lf_\sigma NW_\tau M) = r_S(LW_\tau Nf_\tau M)$, which implies 
$r_S(r_{L\sigma M'}(D)) = r_S(r_{L'\tau M}(D))$. $\Box$

The following corollary may be reserved for the reader to prove.

\begin{corollary}
Let $k<X>$ be a free associative algebra, and "$\leq$" a semigroup partial ordering of $<X>$ with the descending chain condition.

If $S$ is a reduction system on $k<X>$ compatible with $\leq$ and having no ambiguities, then the set of $k$-algebra relations
$W_\sigma =f_\sigma$ ($\sigma \in S$) is independent.

More generally, if $S_1 \subseteq S_2$ are reduction systems, such that $S_2$ is compatible with $\leq$ and all its 
ambiguities are resolvable, and if $S_2$ contains some $\sigma$ such that $W_\sigma$ is irreducible with respect to $S_1$, 
then the inclusion of ideals associated with these systems, $I_1 \subseteq I_2$, is strict. 
\end{corollary}

\section{Remarks}

First I would like to recommend the reader to read the original paper \cite{bergman1} of Bergman, 
because it is written with a broad perspective over a 
lot of algebraic structures, where the reader will find many interesting materials.  

I have to mention at least that there is the correction and updates for the paper. Refer to Bergman \cite{bergman2}.

The Diamond Lemma has another origin, although Newman \cite{newman} is already mentioned.  For that, refer to 
Bokut et al \cite{bob} and Shirshov \cite{shirshov1}. Also see Matveev \cite{Matveev}. Historically, Shirshov \cite{shirshov1} gave the present lemma 
first for Lie algebras. Someone calls the Diamond Lemma Schirshov-Bergman's diamond lemma.

To get more recent trends for the Diamond Lemma, I would like to cite, among others, Chenavier \cite{chena1}, 
Chenavier and Lucas \cite{chena-lucas}, Elias \cite{elias1} and Tsuchioka \cite{tuchioka1}. There the reader will find much more 
information on the lemma and see some practices as an application to representation theory, as an example. 

Recently, there is a trend about Composition-Diamond Lemma. But I do not touch on it here. I shall have an opportunity in another occasion.

I think the Diamond Lemma and its techniques can be applied not only to mathematics but also to many scientific fields.  

\theendnotes \setcounter{endnote}{0}

\bigskip 

%\centerline{\bf \large Appendix \rm} 

%\appendix 

%%%%%%%%%%%%%%%%%%%%%%%%%%%%%%%%%%%%%

\section*{Appendix}

In this appendix, we show a well-known basic theorem for modules without a proof, which is used, in the proof of 
Theorem 3.1, namely the Diamond Lemma, in the third section.

\begin{theorem}
Let $R$ be a ring and $0 \longrightarrow A \stackrel{f}{\longrightarrow} B  \stackrel{g}{\longrightarrow}  C \longrightarrow 0$ a short 
exact sequence of $R$-module homomorphism. Then the following statements are equivalent.

$(1)$ There is an $R$-module homomorphism $h : C \longrightarrow B$ with $gh = id_C$.

$(2)$ There is an $R$-module homomorphism $k : B \longrightarrow A$ with $kf = id_A$.

$(3)$ The given sequence is isomorphic (with identity maps on $A$ and $C$) to the direct sum short exact sequence
$$0 \longrightarrow A \longrightarrow A \oplus C \longrightarrow C \longrightarrow 0,$$
in particular $B \simeq A \oplus C$.
\end{theorem}

\noindent \sc Proof\rm . See for example Hungerford \cite{hungerford}. $\Box$

%%%%%%%%%%%%%%%%%%%%%%%%%%%%%%%%%%%%%

\section*{Ackowledgemant}

I would like to thank Professor Emeritus Dr. Michiel Hazewinkel for correcting a preliminary version  of this note 
and his encouragements to me.

\end{document}